\documentclass[letterpaper, 10 pt, conference]{ieeeconf}
\IEEEoverridecommandlockouts                              % This command is only
 \usepackage{theorem}

\usepackage{enumerate}
 \usepackage{graphics} % for pdf, bitmapped graphics files
\usepackage{epsfig} % for postscript graphics files
\usepackage{url}
\usepackage{amsmath} % assumes amsmath package installed
 \usepackage{amssymb}  % assumes amsmath package installed
\usepackage{verbatim} % for comment
\usepackage{setspace} % for comment

\newcommand{\bea}{\begin{eqnarray}}
\newcommand{\eea}{\end{eqnarray}}
\newcommand{\beas}{\begin{eqnarray*}}
\newcommand{\eeas}{\end{eqnarray*}}
\newcommand{\be }{\begin{equation}}
\newcommand{\ee }{\end{equation}}

\newcommand{\bi}{\begin{itemize}}
\newcommand{\ei}{\end{itemize}}

\newcommand{\dist}{\operatorname{dist}}
\newcommand{\diag}{\operatorname{diag}}

\newcommand{\Int}{\operatorname{{\mathrm Int}}}
\newcommand{\erf}{\operatorname{{\mathrm erf}}}

\newcommand {\R}{\mathbb R}

\newtheorem{Theorem}{Theorem}
\newtheorem{Definition}{Definition}
\newtheorem{Lemma}{Lemma}
\newtheorem{Proposition}{Proposition}
\newtheorem{Corollary}{Corollary}

 \newtheorem{Example}{Example}
\newtheorem{Remark}{Remark}

\newcommand{\sname}{}
\newcommand{\slabel}[1]{\debug{\fbox{\tiny \sname #1}}\label{\sname #1}}
\newcommand{\debug}[1]{}              % final version
\newcommand{\FB}{\begin{figure}[t]\centering}
\newcommand{\FE}[2]{\caption{#2
\debug{\fbox{\sname #1}}} \slabel{#1} \end{figure}}
\newcommand{\tB}{\begin{table}[hbtp]\centering}
\newcommand{\tE}[2]{\caption{#2
\debug{\fbox{\sname #1}}}\slabel{#1} \end{table}}

 \newcommand{\sostfull}{contraction after a small overshoot and short transient}
 \newcommand{\sostshort}{SOST}

 \newcommand{\stfull}{contraction after a  short transient}
\newcommand{\stshort}{ST}

\newcommand{\soshort}{SO}
\newcommand{\sofull}{contraction after a  small overshoot}

\newcommand{\wefull}{weakly expansive}
\newcommand{\sweshort}{WE}

\title{\Large On Three Generalizations of Contraction}

\author{Eduardo D. Sontag  and Michael Margaliot and Tamir Tuller% <-this % stops a space
\thanks{EDS (sontag@math.rutgers.edu) is with the
   Dept. of Mathematics and the Cancer Institute of New Jersey, Rutgers University, Piscataway, NJ 08854, USA;
  MM (michaelm@eng.tau.ac.il) is with the School of Electrical Engineering-Systems
and the Sagol School of Neuroscience, Tel
Aviv University, Israel 69978;
TT (tamirtul@post.tau.ac.il) is with the Dept. of Biomedical Engineering and the Sagol
School of Neuroscience,  Tel-Aviv University, Tel-Aviv 69978, Israel.
EDS's work is supported in part by grants NIH 1R01GM086881 and 1R01GM100473,
and AFOSR FA9550-11-1-0247.
 The research of MM is partly supported by the ISF.
}}

\begin{document}
\maketitle
\thispagestyle{empty}
\pagestyle{empty}

\begin{abstract}
%%%%%%%%%%%%%%%%%%%%%%%%%%%%

%Contraction theory has found numerous applications in
%the analysis and design of non-linear systems.
%Contraction is usually used to prove asymptotic properties such as
%convergence to an attractor  or entrainment to a periodic excitation.

We introduce three  forms of  generalized contraction~(GC).
Roughly speaking, these  are motivated by allowing
   contraction to take place
     after
small transients in time and/or amplitude.
Indeed,  contraction is usually
used to prove asymptotic  properties, like
convergence to an attractor  or entrainment to a periodic excitation,
and
   allowing initial transients does
    not affect this asymptotic behavior.

We provide sufficient conditions
for  GC, and
       demonstrate their usefulness
    using
     examples of
  systems    that
    are not contractive, with respect to any norm, yet are~GC.
%%%
\end{abstract}

%\begin{keywords} Contraction, stability, entrainment, phase locking, systems biology,  ribosome %flow model. \end{keywords}}

\section{Introduction}
%%%%%%%%%%%%%%%%%%%%%%%%%%%%%%%%%%%%
 A  dynamical  system is called contractive
  if  any two   trajectories   converge to one other at an exponential rate. This implies many desirable
  properties including convergence to an attractor (if it exists),
  and entrainment to   periodic excitations~\cite{LOHMILLER1998683,entrain2011}.

Contraction theory is a powerful tool for analyzing nonlinear dynamical systems,  with applications in
control theory~\cite{cont_mech},
observer design~\cite{observer_posi_2011},
synchronization of coupled oscillators~\cite{wang_slotine_2005}, and more.
Recent extensions include: the notion of partial contraction~\cite{partial_cont},
analyzing a network of interacting
 contractive  elements~\cite{russo_hier},
 a Lyapunov-like characterization of incremental stability~\cite{angeli_inc},
 and
a LaSalle-type  principle for contractive systems~\cite{contra_sep}.
%%%
A contractive system with added
diffusion terms or random noise  still satisfies certain asymptotic properties~\cite{Aminzare201331,slotine_cont_noise}.
In this respect,
contraction is  a robust property.

In this paper, we introduce three forms of generalized contraction~(GC).
  These
  are motivated by requiring contraction to take place
   only after
  arbitrarily small
transients in time and/or amplitude. Indeed,   contraction is usually
used to prove \emph{asymptotic} properties, and thus allowing (arbitrarily small)
transients seems reasonable.
 We demonstrate the usefulness of these generalizations
  by showing
 several  examples
 of systems that are \emph{not} contractive with respect to any norm, yet are
 a GC.

\begin{comment}
We use standard  notation.

$A'$ is the transpose of the matrix~$A$.
$\sigma(A)$ denotes the set of eigenvalues of~$A$.
 The \emph{spectral abscissa}
 of~$A \in \R^{n\times n}$ is
\[
\alpha (A):= \max   \{ \mbox{Re}( \lambda  ) : \lambda \in \sigma(A) \},
\]
and the  \emph{spectral radius}
  is
\[
\rho (A):= \max   \{ | \lambda  | : \lambda \in \sigma(A) \}.
\]
  $A $ is called \emph{Hurwitz}  if~$\alpha(A)<0$.
 It is called \emph{Metzler} if all its off-diagonal entries
are non-negative. $A$ is said to be \emph{irreducible}
if for every non-empty proper subset~$I$ of $ N:=\{1,\dots,n\}$
there is an~$i\in I$ and a~$j \in (N \setminus I)$
such that~$a_{ij} \not = 0$.

\end{comment}

The remainder of this paper is organized as follows.
Section~\ref{sec:cthe} provides a brief  overview of contraction theory.
Section~\ref{sec:main} describes our main results. The proofs of these results are detailed
in Section~\ref{sec:proofs}.
Section~\ref{sec:examples} demonstrates
the results
using  a simple model of a
biochemical control system.

%%%%%%%%%%%%%%%%%%%%%%%%%%%%%%%%%%%%%%%%%%%%%%
\section{Contraction Theory} \label{sec:cthe}
%%%%

  %%%%%%%%%%%%%%%%%%%%%%%%%%%%%
 We begin with a brief review of some
 ideas from contraction theory.
 For more details, including the historic development of contraction theory,
 and the relation to other notions, see e.g.~\cite{soderling_survey,cont_anc,RufferWouwMueller:2013:Convergent-Systems-vs.-Incremental-Stabi:}.

 Consider the time-varying
 system
 \be\label{eq:fdyn}
            \dot{x}=f(t,x),
 \ee
evolving on a  convex set~$\Omega \subset \R^n $. We assume that~$f(t,x)$ is differentiable with respect to~$x$, and that
both~$f(t,x)$ and~$J(t,x):=\frac{\partial f}{\partial x}(t,x)$ are continuous in~$(t,x)$.
Let~$x(t,t_0,x_0)$   denote the solution of~\eqref{eq:fdyn}
 at time~$t \geq t_0$ with~$x(t_0)=x_0$
 (for the sake of simplicity, we assume from here on
 that
 $x(t,t_0,x_0)$ exists and is unique for all~$t \geq t_0\geq 0$
 and all~$x_0 \in \Omega$).

 Recall that~\eqref{eq:fdyn} is called
  \emph{contractive}~\cite{LOHMILLER1998683} on~$\Omega$ with respect to  a norm~$|\cdot| :\R^n \to \R_+$
 if there exists~$c>0$
 such that
 \be\label{eq:contdef}
            |x(t_2,t_1,a)-x(t_2,t_1,b)|  \leq   \exp( -   (t_2-t_1) c ) |a-b|
 \ee
 for all $t_2\geq t_1 \geq   0 $ and all~$a,b \in \Omega$.

 In other words,   any two trajectories contract to one another at an exponential rate. This implies
 in particular that the initial condition is ``forgotten''.

Recall that a vector norm~$|\cdot|:\R^n \to \R_+$
 induces a matrix measure
$\mu :\R^{n\times n} \to \R$ defined by
$
            \mu(A):=\lim_{\epsilon \downarrow 0} \frac{1}{\epsilon}
            (||I+\epsilon A||  -1),
$
where~$||\cdot||:\R^{n\times n}\to \R_+$ is the matrix norm induced by~$|\cdot|$.
%%%
A standard  approach for
proving contraction is based on  bounding some matrix measure of
 the Jacobian~$J  $. (This is in fact a particular case of using a Lyapunov-Finsler
 function to prove contraction~\cite{contra_sep}).
\begin{Theorem}~\cite{entrain2011}
If there exists a vector norm~$|\cdot|$ and~$c>0$ such that the induced
 matrix measure~$\mu:\R^{n\times n} \to \R$ satisfies
 \be\label{eq:jtc}
 %%%%
            \mu(J(t,x)) \leq -c,
%%%
\ee
for all $t_2\geq t_1 \geq   0 $ and all~$x \in \Omega$
then~\eqref{eq:contdef}
holds.
\end{Theorem}

One important  implication of
contraction is \emph{entrainment} to a periodic excitation.
%%%%%
Recall that~$f:\R_+ \times  \Omega  \to \R$ is called~\emph{$T$-periodic} if
 \[
 f(t,x) = f(t + T,x)
  \]
  for all $t \geq 0$ and all~$x \in \Omega$.
 %%%%%%
\begin{Theorem}~\cite{entrain2011}
If~\eqref{eq:fdyn} is contractive
and~$f$  is~$T$-periodic
then  there exists  a unique periodic solution $\alpha:[0,\infty]\to \Omega$ of~\eqref{eq:fdyn},
of period~$T$, and
\[
            \lim_{t\to \infty}|x(t,0,a )- \alpha(t) | =0, \text{ for all } a\in  \Omega.
\]
\end{Theorem}

In other words, every trajectory~$x(t,0,a )$ converges to the unique
periodic solution.
Entrainment   is important  in various applications
ranging from biological systems~\cite{RFM_entrain,entrain2011}
to the stability of the power grid~\cite{dorf-bullo}.

The next section presents our main results. All the proofs are placed in Section~\ref{sec:proofs}.

\section{Main Results}\label{sec:main}
%%%%%%%%%%%%%%%%%%%%%%%%%%%%%%%%%%%%%%

We begin by defining  three generalizations of
contraction.
%%%
 \begin{Definition}\label{def:qcont}
 %%%%%
 The time-varying system~\eqref{eq:fdyn} is said to be:
\begin{itemize}
\item  a \emph{\sostfull}~(\sostshort) on~$\Omega$
with respect to  a norm~$|\cdot| :\R^n \to \R_+$ if for each~$\varepsilon >0$
and each~$\tau>0$
  there exists~$\ell=\ell(\tau,\varepsilon)>0$
  %%%%% and~$\varepsilon=\varepsilon(T)>0$, with~$\varepsilon \to 0$ as~$T  \to 0$,
 such that
 \begin{align}\label{eq:qcont}
 %%%%%%
            |x(t_2+\tau,& t_1,a)-   x(t_2+\tau,t_1,b)| \nonumber \\&  \leq   (1+\varepsilon) \exp(-  (t_2-t_1) \ell ) |a-b| \,
 \end{align}
 %%%%%%
for all $t_2\geq t_1\geq 0$  and all $a,b \in \Omega$.
%%%%%%%%%%%%
 \item a  \emph{\sofull} ({\soshort}) on~$\Omega$
with respect to  a norm~$|\cdot| :\R^n \to \R_+$ if for  each~$\varepsilon >0$
  there exists~$\ell=\ell(\varepsilon)>0$
 such that
 \begin{align}\label{eq:uninew}
 %%%
            |x(t_2, & t_1,a)-x(t_2,t_1,b)| \nonumber \\&
             \leq   (1+\varepsilon) \exp(-  (t_2-t_1) \ell ) |a-b|
 \end{align}
 for all  $t_2\geq t_1\geq 0$ and all  $a,b \in \Omega$.
%%%%%%%%%%%%
%%%%%%%%%%%%%
\item    a \emph{\stfull} (\stshort) on~$\Omega$
with respect to  a norm~$|\cdot| :\R^n \to \R_+$ if for
 each~$\tau>0$
  there exists~$\ell=\ell(\tau)>0$
  %%%%% and~$\varepsilon=\varepsilon(T)>0$, with~$\varepsilon \to 0$ as~$T  \to 0$,
 such that
 \begin{align}\label{eq:ucont}
            |x(t_2+\tau, & t_1,a)-x(t_2+\tau,t_1,b)|  \nonumber \\&
            \leq     \exp(-  (t_2-t_1) \ell ) |a-b|
 \end{align}
 %%%
 \end{itemize}
  for all $t_2\geq t_1\geq 0$  and all $a,b \in \Omega$.
 \end{Definition}

It is clear    that every contractive
 system is~\sostshort,
{\soshort}, and~{\stshort}. Thus, all these notions
  are
    generalizations of   contraction.

  The motivation for these definitions stems from the fact that   important
    applications of contraction are  in proving \emph{asymptotic} properties.
    For example,
    proving that an equilibrium point is globally attracting
    or that the state-variables entrain to a  periodic excitation.
    These properties describe what happens as~$t\to \infty$, and so
    it seems natural to generalize contraction in a way that allows
    initial  transients in time and/or amplitude.

In particular, the definition of {\sostshort}  is motivated
by  requiring   contraction  at an exponential rate,
    but only
    after an (arbitrarily small) time~$\tau$, and with  an (arbitrarily small)   overshoot~$(1+\varepsilon)$.
    However, as we will see below when the convergence rate~$ \ell$
    may depend  on~$\varepsilon$
      a somewhat  richer  behavior may occur.
    %%%%%%%%
The definition of {\soshort} is similar to that of~\sostshort, yet now the convergence
rate~$\ell$ depends only on~$\varepsilon$, and there is
no time transient~$\tau$ (i.e.,~$\tau=0$). In other words,~{\soshort} is a uniform (in~$\tau$)
version of~\sostshort.
It is clear that~{\soshort} implies~{\sostshort} and we will see
below that under a mild technical
condition on~\eqref{eq:fdyn} {\soshort}
and   {\sostshort}  are equivalent.

\begin{comment}
%%%%%%%%%%%%%%%%%%%%
\begin{Remark}\label{rem:cont_jaco}
%%
Let~$\alpha(A)$ denote the maximum
over the real parts of all the eigenvalues of~$A$.
It is well-known that if~\eqref{eq:fdyn} is contractive
 on~$\Omega$ with respect to some norm~$|\cdot|$
then~$\alpha(J(t,x))<0$ (i.e.,~$J(t,x)$ is a Hurwitz matrix) for all~$t\geq 0$ and all~$ x \in \Omega$.
It follows from Remark~\ref{re:nexpan} that if~\eqref{eq:fdyn} satisfies anyone of
 the three forms of GC on~$\Omega$
then~$\alpha(J(t,x))\leq 0 $  for all~$t\geq 0$ and all~$x \in \Omega$.
(Indeed, if~$\alpha(J(\bar t,\bar x))> 0 $ for some~$\bar t \geq 0$, $\bar x \in \Omega$
then a linear analysis in the vicinity of this point implies that there exist two trajectories
that diverge from one another
and this contradicts~\eqref{eq:exp}.)
\end{Remark}
 %%%%
\end{comment}

The next example shows that {\stshort} is \emph{not} equivalent to  contraction.
Recall that the \emph{error function} is defined
as~$\erf(z):=\frac {2}{\sqrt{\pi} } \int_0^ z \exp(-s^2) \mathrm{d} s $.
%%%
%%%
\begin{Example}\label{exa:scalarsys}
%%%%
            Consider the scalar time-varying  system
            \be\label{eq:scals}
            %%%%
                    \dot x(t)=(\exp(-t^2)-1)x(t)
            \ee
            evolving on~$\Omega:=(-1,1)$.
           It is straightforward to
           show that this system is  {not}   contractive
           with respect to \emph{any} norm
           (note that the Jacobian~$J(t)=\exp(-t^2)-1$  satisfies~$J(0)=0$).
           Yet,~\eqref{eq:scals}
            is~{\stshort}. To show this, pick~$\tau>0$.
           Note that
  for all~$t\geq t_1$, $x(t,t_1,a )= f(t,t_1) a$,
            where
            \be\label{eq:explisol}
            f(t,t_1):=\exp \left(g(t)-g(t_1)\right)   ,
              \ee
 %%%
with
            \be\label{eq:defgt}
            g(t):= \frac{\sqrt{\pi}}{2}  \erf(t)-t .
            \ee
  Thus, for any norm~$|\cdot|:\R\to \R_+$,
            \begin{align*}
            |x(t_2,t_1,a)-x(t_2,t_1,b)| &=f(t_2,t_1)  |a-b|,
            \end{align*}
%%%
  We need to show that
    there exists~$\ell_1=\ell_1(\tau   )>0$ such that
\[
   f(t_2+\tau,t_1 ) \leq   \exp(-(t_2-t_1) \ell_1 ),\quad \text{for all } t_2\geq t_1\geq 0,
\]
 or, equivalently, that
\be\label{eq:whintk}
   a(z,t_1)        \leq 0, \quad \text{for all } t_1,z \geq 0,
\ee
where $z:=t_2-t_1 $, and~$a(z,t_1):=g(t_1+\tau+ z) -g(t_1) +  z \ell_1 $.
Since
 \be\label{eq:gder}
                        \dot{g}(t)= \exp(-t^2)-1 ,
\ee
  we have
 \begin{align*}
                  \frac{d}{dz} a(z,t_1)&=  \exp( - (t_1+\tau+ z )^2 )  -1 +\ell_1 \\
                  &\leq   \exp( - \tau^2 ) -1  +\ell_1.
\end{align*}
 Taking
$
 \ell_1 :=(1-\exp( - \tau^2 ))/2>0$,
 %%%%
  yields~$\frac{d}{dz} a(z,t_1)<0$ for all~$t_1,z \geq 0$, so \begin{align*} a(z,t_1)&\leq a(0,t_1)\\
       &= g(t_1+\tau ) -g(t_1)  \\
       &\leq 0       ,
\end{align*} where the last inequality follows from~\eqref{eq:gder}.
We conclude that~\eqref{eq:whintk} indeed holds, so~\eqref{eq:scals} is \stshort.
%%%%%%
\end{Example}

For  \emph{time-invariant}  systems evolving on  a compact set
 it is possible to give a simple  sufficient condition for~{\stshort}.
Let~$\Int(S)$ denote the interior of a set~$S$.
We require the following definitions.
%%%%%
\begin{Definition}
We say that~\eqref{eq:fdyn} is
\emph{non expansive}~(NE) with respect to the norm~$|\cdot|$ if
for all~$a,b\in\Omega$ and all~$s_2 > s_1\geq 0$
   \be\label{eq:exp}
        |x(s_2,s_1,a)-x(s_2,s_1,b)|   \leq    |a-b|.
\ee
%%%%%%
We say that~\eqref{eq:fdyn} is \emph{weakly contractive}~(WC) if~\eqref{eq:exp}
holds with~$\leq $ replaced by~$<$.
%%%%%%%%
\end{Definition}

\begin{Definition}
%%%
The  time-invariant system
\be\label{eq:time_in_var}
%%%%%%%%%%%%%%%%%%%%%%%%%%%%%%%%%%%%%%%
                    \dot{x}=f(x),
\ee
 evolving on a   compact and convex set~$\Omega \subset \R^n$,
 is said to be \emph{interior contractive}~(IC)
 if it satisfies the
following properties:
\begin{enumerate}[(a)]
\item   for every~$x_0 \in \partial \Omega$,
    $x(t,x_0) \not \in \partial \Omega$ for all~$t>0$.
%%%%%
\item  there exists a matrix measure~$\mu:\R^{n\times n}\to \R$ such that
\be \label{eq:mucom_inv}
\mu(J(x))<0,\quad \text{for all } x \in \Int(\Omega).
\ee
\end{enumerate}
\end{Definition}
%%%%%%%

Note that conditions~(a) and~(b)  do not necessarily imply
contraction on~$\Omega$, as it is possible that~$\mu(J(x))=0$   for some~$x\in \partial \Omega$.
Yet,~\eqref{eq:mucom_inv} does   imply that~\eqref{eq:time_in_var} is
non-expansive on~$\Omega$.

\begin{Theorem}\label{thm:time_invar}
%%%
If the system~\eqref{eq:time_in_var} is IC then
 it is {\stshort}.
\end{Theorem}

As noted above, the introduction of the GC forms
  is motivated by the idea that contraction is used to prove asymptotic results,
so allowing
initial transients should increase the class of systems that can be analyzed while still
allowing to prove asymptotic results.
The next result demonstrates this.
%%%%%%%%%%%%%%%%%%%%%%%%%%%%%%%%%%%%%%%%%%%
\begin{Corollary}\label{coro:attract}
      Suppose that~\eqref{eq:time_in_var} is~IC.
  Then~\eqref{eq:time_in_var} admits
  a unique equilibrium point~$e \in \Int(\Omega)$,
and~$\lim_{t \to \infty}x(t,a)= e $ for all~$a\in \Omega$.
\end{Corollary}

One may perhaps expect that we can generalize  Theorem~\ref{thm:time_invar}
to the time-varying case as well, that is,
that if  the
time-varying system~\eqref{eq:fdyn},
 evolving in a   compact and convex set~$\Omega \subset \R^n$,
 satisfies:
\begin{enumerate}[(a)]
\item   for every~$x_0 \in \partial \Omega$ and every~$t_1 \geq 0$,
    \be\label{eq:prop1_co}
    x(t,t_1,x_0) \not \in \partial \Omega, \text{ for all } t>t_1,
    \ee
%%%%%
\item  there exists a matrix measure~$\mu:\R^{n\times n}\to \R $ such that
    \be\label{eq:prop2_co}
\mu(J(t,x))<0, \text{ for all } x \in \Int(\Omega), \text{ and all } t\geq t_1\geq 0,
\ee
\end{enumerate}
then~\eqref{eq:fdyn} is {\stshort} on~$\Omega$.
However, the next example shows that this is not so.
\begin{Example}\label{exa:xounter}
Consider the  {scalar} system
\[
            \dot{x}(t)=-\frac{x(t)}{t+1}   ,
\]
evolving in~$\Omega:=[-1,1]$. The Jacobian   is~$J(t,x)=-(t+1)^{-1}$,
and   properties~\eqref{eq:prop1_co} and~\eqref{eq:prop2_co} hold.
Yet, this system is \emph{not} {\sostshort} on~$\Omega$ (and, therefore,
it is clearly not~{\stshort} on~$\Omega$). Indeed, assume otherwise. Pick~$\tau,\varepsilon>0$.
Then there exists~$\ell=\ell(\tau,\varepsilon)>0$ such that~\eqref{eq:qcont}
holds. Since~$x(t)=(t+1)^{-1} x(0)$, Eq.~\eqref{eq:qcont} with the particular choice~$t_1=0$
implies that
\[
   (t_2+\tau+ 1)^{-1}    \leq (1+\varepsilon)\exp(-\ell t_2)  ,\quad \text{for all } t_2\geq 0,
\]
i.e.,
\[
       \exp(\ell t_2) \leq     (1+\varepsilon) (t_2+\tau+ 1)     , \quad \text{for all }t_2 \geq 0,
\]
but this clearly  cannot hold for~$t_2>0$ sufficiently large.
\end{Example}

To provide a sufficient
 condition for generalized contraction
 of
    the \emph{time-varying} system~\eqref{eq:fdyn},
  we require the following definition.
  %%%%
 \begin{Definition}
 %%%%%
 System~\eqref{eq:fdyn} is said
  to be \emph{nested contractive}~(NC) on~$\Omega$ with respect to a norm~$|\cdot|$
if there exist  convex   sets~$\Omega_\zeta \subseteq \Omega$,    and norms~$|\cdot|_\zeta:\R^n\to \R_+$,
 where~$\zeta\in (0,1/2]$,
 such that the following conditions hold.
 \begin{itemize}
                    \item $\cup_{\zeta \in (0,1/2]} \Omega_\zeta=\Omega$, and
                    \be\label{eq:setsinc}
                    \Omega_{\zeta_1} \subseteq \Omega_{\zeta_2},\quad \text{for all } \zeta_1 \geq \zeta_2.
                    \ee
                    \item For every~$ \tau>0$ there exists~$\zeta=\zeta( \tau)\in(0,1/2]$,
                    with~$\zeta ( \tau)\to 0$ as~$ \tau \to 0 $,
                    such that  for every~$a\in \Omega$ and every~$t_1 \geq 0$
                    \be \label{eq:enter}
                            x(t ,t_1,a)\in\Omega_{\zeta },\quad \text{for all } t\geq t_1+\tau,
                     \ee
                     and~\eqref{eq:fdyn}
                                              is contractive  on~$\Omega_\zeta$ with respect to~$|\cdot|_\zeta$.
                     %%%
                     \item The  norms  $|\cdot|_\zeta$ converge to~$|\cdot| $ as~$\zeta\to 0$, i.e., for every~$\zeta>0$ there exists~$s =s(\zeta) >0$, with~$s(\zeta) \to 0$
                          as~$ \zeta \to 0$, such that
                         \[
                                  (1-s)|  y | \leq |y|_\zeta  \leq (1+s) |y|      ,\quad \text{for all } y \in \Omega       .
                         \]
                      \item System~\eqref{eq:fdyn}
                                              is  non-expanding with respect
                     to~$|\cdot|$ on~$\Omega$.

 \end{itemize}
  %%%%%
 \end{Definition}

Eq.~\eqref{eq:enter} means that
 after an arbitrarily  short time every trajectory enters and remains in
 a subset~$\Omega_\zeta$ of the state space on which we have contraction
 with respect to~$|\cdot|_\zeta$.
%%%

 \begin{Theorem}\label{thm:qcon}
 %%%%%
 If the system~\eqref{eq:fdyn} is NC then it is
  {\sostshort}.
 %%%%%
 \end{Theorem}
 %%%

The next example demonstrates the usefulness  of Theorem~\ref{thm:qcon}
by using it to prove that the system in Example~\ref{exa:scalarsys} is~{\sostshort}
    \emph{without}
 using the explicit solution~\eqref{eq:explisol}.
%%%% *****************************************************
\begin{Example}
%%%%
            Consider again the scalar time-varying  system~\eqref{eq:scals}.
            Fix arbitrary~$t_1\geq 0$, $a_1\in (-1,1)$, and  rewrite~\eqref{eq:scals}   as
            \begin{align}\label{eq:expsys}
                    \dot x_1 &=(\exp(-x_2^2)-1)x_1 , & x_1(t_1)=a_1,\nonumber\\
                    \dot x_2 &=1, &x_2(t_1)=t_1,
            \end{align}
evolving on~$\Omega :=\{ x\in \R^2:x_1 \in (-1,1),  x_2 \geq t_1 \}$.
 Note that any two feasible initial conditions~$a,b \in \Omega $ for this systems
  satisfy $a_2=b_2=t_1$.
%%%%
The Jacobian of~\eqref{eq:expsys} is \[
            J(x)=\begin{bmatrix}    -1+\exp(-x_2^2) &-2x_1x_2\exp(-x_2^2)\\
                                      0              & 0
             \end{bmatrix}.
\]
For any~$\zeta \in (0,1/2]$, let \[ \Omega_\zeta := \{x\in \Omega : x_1\in(-1,1), \;x_2\geq t_1+\zeta  \}, \] and let~$|\cdot|_\zeta:=|\cdot|_1$, that is, the~$L_1$ norm.
Note that~\eqref{eq:enter} holds with~$\zeta(\tau):=\min\{\tau,1/2\}$, and
  that for every~$x \in \Omega_\zeta$,
\begin{align}\label{eq:j11}
                    J_{11}(x) & =  -1+\exp(-x_2^2)\nonumber\\
                              &\leq -1+\exp(-(t_1+\zeta)^2)\nonumber\\
                               &< -1+\exp(- \zeta ^2)\nonumber\\
                               &<0.
\end{align} Let~$d_i(t_2,t_1,a,b):=|x_i (t_2,t_1,a)-x_i(t_2,t_1,b)|$, $i=1,2$. Then \begin{align*}
                    |x(t_2,t_1,a)-x(t_2,t_1,b)|_1 &=d_1(t_2,t_1,a,b)+d_2(t_2,t_1,a,b)\\
                                &=d_1(t_2,t_1,a,b)+ |t_2- t_2|\\
                                &=|x_1 (t_2,t_1,a)-x_1(t_2,t_1,b)|.
\end{align*}
 Combining this with~\eqref{eq:j11} implies that all the conditions in Theorem~\ref{thm:qcon} hold, so we conclude that~\eqref{eq:expsys} is {\sostshort}  with respect to
the~$L_1$ norm.
\end{Example}

The next section describes, using a specific
 mathematical model,
 one possible mechanism for~{\stshort}. Namely,
as we change the parameters in a contractive system, it may become~{\stshort}
when it hits the ``verge'' of contraction.
For two vectors~$a,b\in\R^n$, we write~$a\geq b$   if~$a_i\geq b_i$
for~$i=1,\ldots,n$. A matrix~$M \in\R^{n\times n}$ is called \emph{Metzler} if~$m_{ij}\geq 0$ for all~$i \not = j$.

%%%%%%%%%%%%%%%%%%%%%%%%%%%%%%%%%%%%%%%%%%%%
%%%%%%%%%%%%%%%%%%%%%%%%%%%%%%%%%%%%%%%%%%%%%
\section{An application: A biochemical control~circuit}\label{sec:examples}
Consider the system
\begin{align}\label{eq:bio}
\dot{x}_1&=g(x_n)-\alpha_1 x_1,\nonumber\\
\dot{x}_2&=x_1-\alpha_2 x_2\nonumber,\\
\dot{x}_3&=x_2-\alpha_3 x_3 \nonumber,\\
       &\vdots\nonumber\\
\dot{x}_n&=x_{n-1}-\alpha_n x_n  ,
\end{align}
where~$\alpha_i>0$, and
\[
g(u):=\frac{1+u}{k+u} , \quad \text{with } k>1.
\]
As explained  in~\cite[Ch.~4]{hlsmith} this may model a simple
  biochemical control circuit for protein synthesis in the cell.
%%
%%%
The~$x_i$s represent concentrations of various macro-molecules in the cell
and therefore must be non-negative. It is straightforward to
verify that~$x(0) \in \R^n_+$ implies that~$x(t) \in \R^n_+$ for all~$t \geq 0$.
%%%%%%%%%%%%%%%%%%%%%
%%%%%%%%%%%%%%%%%%%%%%
\begin{Proposition}\label{prop:bio}
%%%
Let~$\alpha:=\prod_{i=1}^n \alpha_i  $.
 If
 \be\label{eq:km1}
   {k-1}  < \alpha {k^2}
   \ee
then~\eqref{eq:bio} is a contraction on~$\R_+^n$.
If
$ {k-1} =\alpha {k^2}
 $  then~\eqref{eq:bio} is not a contraction, with respect to any norm,
on~$\R_+^n$, yet it is {\soshort} on~$\R_+^n$.
%%%
\end{Proposition}
%%%%%%%%%%%%%%

Note that
 for all~$x\in\R^n_+$,
\be\label{eq:gderi}
        g'(x_n)= \frac{k-1}{(k+x_n)^2} \leq \frac{k-1}{k^2} =g'(0).
\ee
Thus~\eqref{eq:km1} implies that contraction holds
if and only
if the
``total dissipation'' $\alpha$
is strictly larger   than~$g'(0)$.

Using the fact that~$g(u)<1$ for all~$u\geq 0$ it is straightforward to show that
the set
\[
        \Omega_r:= r  ( [0,\alpha_1^{-1}]\times[0, (\alpha_1 \alpha_2)^{-1}] \times\dots\times [0, \alpha^{-1}]  )
\]
is an invariant set of the dynamics for all~$r\geq 1$.
Combining this with Prop.~\ref{prop:bio}
     implies that~\eqref{eq:bio},
with~$   {k-1}  \leq \alpha {k^2}$,
admits a unique equilibrium point~$e \in \Omega_1 $ and that
\[
            \lim_{t\to \infty} x(t,a)=e,\quad \text{for all } a \in \R^n_+.
\]
This property also follows from a more general result~\cite[Prop.~4.2.1]{hlsmith}
that is proved using the theory of irreducible  cooperative  dynamical systems.
Yet the contraction approach leads to new insights.
For example,
it implies that the distance between trajectories can only decrease,
and can also be used to prove entrainment to   suitable generalizations
of~\eqref{eq:bio} that include periodically-varying inputs.

In the next section, we describe several more related notions and explore the relations between
them.

\section{Additional Notions and Relations}\label{sec:relations}
%%%%%%%%%%%%%%

It is straightforward to show that
each of the three generalizations of  contraction in Definition~\ref{def:qcont}
  implies  that~\eqref{eq:fdyn} is~NE.
 %%%%%%
 One may perhaps expect that any of the three generalizations
 of  contraction in Definition~\ref{def:qcont}
also implies~WC.
By taking~$t_1=s_1$, $\tau=(s_2-s_1)/2>0 $, and~$t_2=s_1+\tau$ in~\eqref{eq:ucont}
it follows   that~{\stshort} does imply~WC.
However, the next example shows that~{\soshort} does not imply WC.
 \begin{Example}\label{exa:eps}
%%%
 Consider the  {scalar}  time-varying system
 %%%%%%%%%%%%%%
            \be\label{eq:shift}
                    \dot{x}(t)=\begin{cases}  0, & 0\leq t \leq 1,\\
                                               -2x(t), &1<t ,
                               \end{cases}
            \ee
            evolving on~$\Omega:=(-1,1)$.
%%%
Clearly, the trajectories of this system are not contracting  for~$t \in [0,1]$.
Yet,  we claim  that  this system is~{\soshort}. To show this,
 pick~$ \varepsilon>0$.
 Let
  \be\label{eq:ellvareps}
                    \ell :=\min\{ \log(1+\varepsilon),1\} .
 \ee
 Note that~$\ell=\ell(\varepsilon) >0$.
To show that~\eqref{eq:uninew} holds, we consider two cases.

\noindent {\sl Case 1.} Suppose that~$t_1 \in [0,1]$.
 In this case,
  the solution of~\eqref{eq:shift} is
  %%%
\be\label{eq:solution}
                                x(t,t_1,a)= \begin{cases}
                                               a, &  t_1 \leq t \leq 1,\\
                                                 \exp(-2(t-1)) a,    &1 \leq t  .
                               \end{cases}
%%%
\ee
Thus,
%%%
\begin{align*}
%%%
                      d:&= |x( t_2  ,t_1, a)-x( t_2  ,t_1 ,b)| \\& = \begin{cases}
                                                |a-b| , &  t_1 \leq t_2   \leq 1,\\
                                                \exp(-2 (t_2 -1 )) | a-b | ,   &1 \leq t_2  .
                               \end{cases}
%%%
%%%
\end{align*}
Let~$r:=(1+\varepsilon)\exp(-\ell(t_2-t_1))|a-b|$.
It follows from~\eqref{eq:ellvareps} that
\begin{align*}
r  & \geq
  (1+\varepsilon)^{1-t_2+t_1}|a-b|,
\end{align*}
so if~$ t_2 \leq 1   $ then clearly~$ d  \leq r $.
%%%

Now suppose that~$ t_2  >  1  $.
If~$1+\varepsilon \geq \exp(1)$ then it follows from~\eqref{eq:ellvareps} that~$\ell=1 $,
so
\begin{align*}
                 r&  =  (1+\varepsilon)\exp(- (t_2-t_1))|a-b| \\
                  &\geq  \exp(1)  \exp(- t_2)|a-b|\\
                  &\geq \exp(2) \exp(-2 t_2 ) | a-b |
                  \\ &=d.
\end{align*}
If~$1+\varepsilon < \exp(1)$ then it follows from~\eqref{eq:ellvareps} that~$\ell=\log(1+\varepsilon) $,
so
\begin{align*}
                 r&  =   (1+\varepsilon)^{1-t_2+t_1}|a-b| \\
                  &\geq   (1+\varepsilon)^{1-t_2 } |a-b|\\
                        &\geq   \exp(1-t_2)  |a-b|\\
                  &\geq \exp(2) \exp(-2 t_2 ) | a-b |
                  \\ &=d.
\end{align*}
Summarizing, in Case~1 we always have~$d \leq r$.

\noindent {\sl Case 2.} Suppose that~$t_1 >1$.
 In this case,
  the solution of~\eqref{eq:shift} is
  %%%
$
                                x(t,t_1,a)=
                                                 \exp(-2(t-t_1)) a,
%%%
$
so
%%%
$
%%%
                      d =          \exp(-2 (t_2 -t_1 )) | a-b |.
%%%
%%%
$
Since~$\ell \leq 1$,~$d\leq  (1+\varepsilon) \exp(-\ell (t_2-t_1)) |a-b| =r$.
Thus, in Case~2 we also have~$d \leq r$, and this proves~{\soshort}.
%%%%%%%%%%%%%%%%%%%%%%%%%%%%%%%%%%%%%%%%%%%%
 %%%%%%%%
\end{Example}

%%%%%%%%%%%%%%%%%%%%%%%%

Summarizing, \eqref{eq:shift}
  is {\soshort} although its trajectories do not
contract  for~$t \in [0,1]$. Clearly, for every fixed~$T>0$
we can build
a system that is  {\soshort} although its trajectories do not
contract  for~$t \in [0,T]$.

\subsection{\sostfull}
%%%%%%%%%%%%%
The next result presents two conditions that are equivalent to~{\sostshort}.
\begin{Lemma}\label{lem:eqdef}
%%%%%
The following conditions are equivalent.
 \begin{enumerate}
 \item \label{item1}
 System~\eqref{eq:fdyn}  is {\sostshort} on~$\Omega$
  with respect to  some vector norm $|\cdot|_v:\R^n\to \R_+  $.
%%%%
\item \label{item2}
 For each~$ {\tau}>0$
  there exists~$ {\ell}= {\ell}( {\tau} )>0$
 such that
 \begin{align}\label{eq:qcontop}
            |x(t_2+\tau,& t_1,a)-  x(t_2+\tau,t_1,b)|_v \nonumber \\
                    &\leq   (1+ {\tau}) \exp(-  (t_2-t_1)  {\ell} ) |a-b|_v ,
 \end{align}
 for all~$t_2\geq t_1\geq 0$ and all~$a,b \in \Omega$.
%%%%%
%%%%
\item \label{item3}
For each~$\varepsilon >0$ and each~$\tau>0$
  there exists~$\ell_1=\ell_1(\tau,\varepsilon)>0$
 such that
 \be\label{eq:qcontnew}
 %%%%%
            |x(t  ,t_1,a)-x(t  ,t_1,b)|_v  \leq   (1+\varepsilon) \exp(-  (t -t_1) \ell_1 ) |a-b|_v ,
 \ee
 for all~$t  \geq t_1+\tau \geq \tau$ and all~$a,b \in \Omega$.
%%%
\end{enumerate}
\end{Lemma}
 %%%%%%%%%%%%%%%%%%%%%%%%%%

 %%%%%%%%%%%%%%%%%%%%%%%%%%%%%%%%%%%%%%%%%%%%
\subsection{\sofull}
%%%%%%%%%%%%%%%%%%

A natural question is under what conditions~{\soshort} and~{\sostshort} are equivalent. To address this issue,
we introduce the following definition.

\begin{Definition}\label{eq:defntr}
%%%
                            We say that~\eqref{eq:fdyn} is \emph{\wefull} (\sweshort)
                            if for each~$\delta>0$ there exists~$\tau_0>0$ such that for all~$a,b \in \Omega$ and all~$t_0\geq 0$
                            \be\label{eq:sep}
                                        | x(t ,t_0,a)-x(t ,t_0,b) |\leq(1+\delta)|a-b|,
                            \ee
 for all $t \in [t_0,t_0+\tau_0]$.
 \end{Definition}

 \begin{Proposition}\label{prop:sepimp}
            Suppose that \eqref{eq:fdyn} is   {\sweshort}.
            Then \eqref{eq:fdyn} is {\sostshort}  if and only if it is {\soshort}.
\end{Proposition}

\begin{Remark}\label{rem:Global_Lip}
Suppose that~$f$ in~\eqref{eq:fdyn} is Lipschitz globally in~$\Omega$
 uniformly in~$t$, i.e. there exists~$L>0$ such that
\[
               |  f(t,x)-f(t,y)|\leq L|x-y|, \quad \text{for all }x,y \in \Omega, \; t\geq t_0.
\]
Then by Gronwall's Lemma (see, e.g.~\cite[Appendix~C]{sontag_textbook})
 \[
        | x(t ,t_0,a)-x(t ,t_0,b) |\leq \exp(L (t-t_0) )|a-b|,
 \]
%%%
  for all $t \geq t_0$, and this implies that~\eqref{eq:sep} holds for~$\tau_0:=\frac{1}{L} \log(1+\delta)>0$.
In particular, if~$\Omega$ is compact and~$f$ is periodic in~$t$ then~{\sweshort} holds under rather weak continuity arguments on~$f$.
\end{Remark}

\begin{figure*}[t]
 \begin{center}
 \input{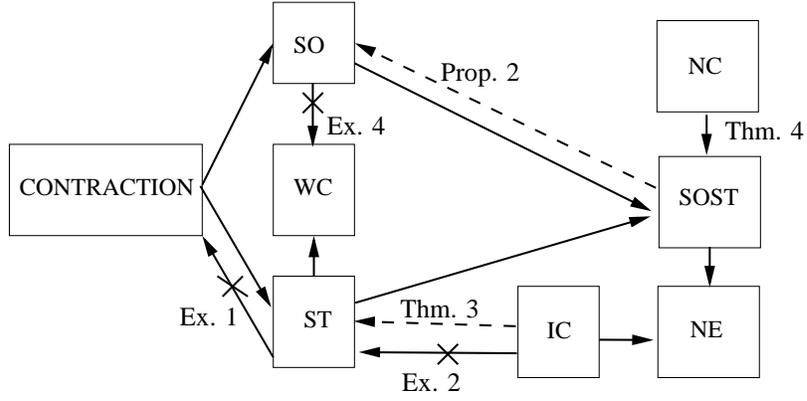}
\caption{Relations between various contraction notions.
An arrow denotes implication;
a crossed out arrow denotes that the implication is in general false;
and a dashed arrow denotes
 an implication that holds under an additional condition.
Some of the relations are immediate.
Others follow from the results   marked  near the arrows.
}
\label{fig:graphbn2}
\end{center}
\end{figure*}

Fig.~\ref{fig:graphbn2} summarizes
the relations between the various contraction
notions.

 \section{Proofs} \label{sec:proofs}
%%%%%%%%%%%%%%%%%%%%%%%%%%%%%%%%%%%%

%%%%%%%%%%%%%%%%%%%%%%%%%%%%%%%%%%%%%%%%%%%%
\noindent {\sl Proof  of Theorem~\ref{thm:time_invar}.}\\
%%%%%%%%%%%%%%%%%%%%%%%%%%%%%%%%%%%%%%%%%%%%
We require the following result.
\begin{Lemma}\label{lem:strict}
%%%%%%%
If   system~\eqref{eq:time_in_var} is IC
then  for each $\tau > 0$ there is a $d > 0$
 such  that
 \[
 \dist(x(t,x_0),\partial \Omega) \geq d ,
 \]
 for all~$x_0 \in \Omega$ and all~$t \geq \tau$.
\end{Lemma}

 {\sl Proof of Lemma~\ref{lem:strict}.}
 %%%%
 Pick  $\tau>0$ and~$ x_0 \in \Omega$.
 The assumption that~$\Omega$ is invariant
 implies that~$\Int(\Omega)$ is also an invariant set of~\eqref{eq:time_in_var}
 (see, e.g.,~\cite[Lemma III.6]{mcs_angeli_2003}).
Combining this with~(a)   implies that
    $x(t,x_0) \not \in \partial \Omega$ for all~$x_0 \in \Omega$ and all~$t>0$, so
%%%%
$e_{x_0}:=\dist(x(\tau,x_0),\partial \Omega)  > 0$.
Thus, there exists a
  neighborhood~$U_{x_0}$ of $x_0$, such that
$\dist(x(\tau,y),\partial \Omega) \geq e_{x_0}/2$
for all $y \in U_{x_0}$.
Cover~$\Omega$ by such~$U_{x_0}$ sets.
 By compactness of~$\Omega$, we can pick a    finite subcover.
  Pick smallest~$e$ in this subcover, and denote
 this by~$d$. Then~$d>0$ and
we have that
$\dist(x(\tau,x_0),\partial \Omega) \geq d$ for all $x_0 \in \Omega$.
Now, pick~$t \geq \tau$ and $x_0 \in \Omega$.
Let $x_1: = x(t-\tau,x_0)$.  Then:
\[
\dist(x(t,x_0),\partial  \Omega) = \dist(x(\tau,x_1),\partial \Omega) \geq d,
\]
and this completes the proof of Lemma~\ref{lem:strict}.~$\square$

To prove  Theorem~\ref{thm:time_invar}, pick~$\tau>0$.
Let~$S_\tau:=\{ x(t,x_0):\; t \geq \tau,x_0 \in \Omega  \}$.
Lemma~\ref{lem:strict} implies that there exists a closed and convex set~$D$
such that
\[
             S_\tau  \subseteq D \subset  \Int(\Omega).
\]
(Note that since~$\Omega$ is convex so is~$\Int(\Omega)$).
 Let~$c_\tau:=\min_{x \in D} \mu(J(x))$. Then~$c_\tau<0$. Thus, the system
 is contractive on~$D $, and
for all~$a,b \in \Omega$ and all~$t\geq 0$
\[
            |x(\tau+t,a)-x(\tau+t,b)|\leq \exp( c_\tau t) |a-b|,
\]
 where~$|\cdot| $
is the vector norm corresponding to the matrix measure~$\mu $.
%%%%
  This establishes~\stshort, and thus
  completes the proof of Theorem~\ref{thm:time_invar}.~$\square$

%%%%%%%%%%%%%%%%%%%%%%%%%%%%%%%%%%%%%%%%%%%%%%%%%

\noindent {\sl Proof of Corollary~\ref{coro:attract}.}\\
%%%%%%%
\noindent Since~$\Omega$ is convex, compact, and invariant, it includes an
 equilibrium point~$e$  of~\eqref{eq:time_in_var}.
By Theorem~\ref{thm:time_invar},
 the system is~{\stshort}. Pick~$a \in \Omega$ and~$\tau>0$.
 Applying~\eqref{eq:ucont}
with~$b=e$ yields
\begin{align*}
            |x(t_2+\tau, & t_1,a)-e|
            \leq     \exp(-  (t_2-t_1) \ell ) |a-e|,
 \end{align*}
 for all~$t_2\geq t_1\geq 0$, where~$\ell>0$. This  completes the proof.~$\square$

%%%%%%%%%%%%%%%%%%%%%%%%%%%%%%%%%%%%%%%%%%%%%%%%%%%
\noindent  {\sl Proof of Theorem~\ref{thm:qcon}.}\\
 %%%%%%%%%%%%%%%%%%%%%%%%%%%%%%%%%%%%%%%%%
Fix arbitrary $\varepsilon>0$ and $t_1\geq 0$.
The function $\zeta=\zeta( \tau) \in (0,1/2]$ is as in the statement of the
Theorem.
For each $ \tau>0$, let $c_{\zeta}>0$ be a contraction constant
on~$\Omega_\zeta$, where we write $\zeta=\zeta( \tau)$ here and in what follows.
Take any $a,b \in \Omega$.
By~\eqref{eq:enter},
$x(t,t_1,a),x(t,t_1,b) \in \Omega_\zeta$ for  all~$t\geq t_1+\tau$, so
 \begin{align*}
|x(t ,t_1,    a)- x(t ,t_1,&b)|_\zeta
 \leq  \exp(-c_{\zeta } (t- t_1-\tau))\\ &\times | x(t_1+\tau ,t_1, a)- x(t_1+\tau ,t_1,b)|_\zeta
 \end{align*}
for all~$t\geq t_1+\tau$.
From the convergence property of norms in the Theorem statement,
there exist~$v_{\zeta },w_{\zeta } >0 $ such
 that
 \begin{align}
                            |y|  & \leq v_{\zeta } | y|_{\zeta } \leq w_{\zeta } v_{\zeta } | y|,\quad \text{for all } y\in\Omega,
 \end{align}
and $v_{\zeta }\rightarrow1$, $w_{\zeta }\rightarrow1$ as $\tau\rightarrow0$.
For~$t\geq t_1+\tau$ let~$p:=t- t_1-\tau$. Then
  \begin{align*}
     |  x(t , & t_1,  a) - x(t ,t_1,b)| \\
          & \leq  v_\zeta  \exp(-c_\zeta p) | x(t_1+\tau ,t_1, a)- x(t_1+\tau ,t_1,b)|_\zeta \\
          &\leq v_\zeta   w_\zeta  \exp(-c_\zeta p) | x(t_1+\tau ,t_1, a)- x(t_1+\tau ,t_1,b)|\\
          &\leq v_\zeta   w_\zeta \exp(-c_\zeta p)    | a-b| ,
  \end{align*}
  where the last inequality follows from the fact
  that   the system is non-expanding with respect to~$|\cdot|$.
Since
$v_{\zeta }\rightarrow1$, $w_{\zeta }\rightarrow1$ as $\tau\rightarrow0$,
$v_{\zeta } w_{\zeta }  \leq  {1+ \varepsilon}$
for $\tau>0$ small enough.
Summarizing, there exists~$\tau_m=\tau_m(\varepsilon)>0$ such that for all~$\tau \in[0, \tau_m]$
  \begin{align}\label{eq:summa}
    |& x(t+\tau ,t_1,  a)- x(t+\tau ,t_1,b)| \nonumber \\
          &\leq  (1+\varepsilon)   \exp(-c_\zeta (t- t_1 ))    | a-b| ,
  \end{align}
for all~$a,b\in \Omega$ and all~$t\geq t_1 $. Now pick~$\tau > \tau_m$.
For any~$t\geq t_1$,
let~$s:= t+\tau-\tau_m $.
Then
  \begin{align*}
    |  x(t+\tau ,t_1,&  a)-   x(t+\tau ,t_1,b)|\\
    &=|  x(s+\tau_m ,t_1,  a)- x(s+\tau_m ,t_1,b)|  \\
    &\leq  (1+\varepsilon)   \exp(-c_\zeta (s- t_1 ))    | a-b| \\
     &\leq  (1+\varepsilon)   \exp(-c_\zeta (t- t_1 ))    | a-b|,
  \end{align*}
and this  completes the proof.~$\square$
%%%%%%%%%%%%%%%%%%%%%%%%%%%%%%%%%%%%%%%%%%%%

\noindent {\sl Proof of Proposition~\ref{prop:bio}.}\\
%%%%%%%%%%%%%%
We require the following result from~\cite{RFM_entrain}.
 \begin{Lemma}\label{lem:per}
%%%%%%%%%%%
Consider a time-varying system
\be \label{eq:gensys}
\dot x = f(t,x)
\ee
evolving on a subset of $X := I_1\times I_2\times \ldots \times I_n\subseteq \R^n_{+}$, where each
$I_j$ is an interval of the form $[0,a]$, $a>0$, or $[0,\infty )$.  Suppose
that the time-dependent vector field $f=\begin{bmatrix}
f_1,\ldots ,f_n\end{bmatrix}'$ has the following
\emph{ boundary-repelling} property:
\bi
\item[{\bf (BR)}]
For each $\delta >0$ and each sufficiently small $\Delta >0$, there exists
$K=K(\delta ,\Delta )>0$ such that, for each $k=1,\ldots ,n$ and each $t\geq 0$,
the condition
\be\label{eq:fk}
%%%%%
x_k\leq \Delta   \text{ and } x_i\geq \delta  ,\quad \text{for every }  1\leq  i \leq  k-1
%%%%
\ee
(for $k=1$, the condition is simply $x_1\leq \Delta $)\\
implies that
\be\label{eq:fkk}
%%%
f_k(t,x) \geq  K,\quad \text{for all } t\geq 0.
\ee
%%%
\ei
Then  given any $\tau>0$ there exists   $\varepsilon  = \varepsilon (\tau)>0$,
 with~$\varepsilon(\tau) \to 0$ as~$\tau \to 0$, such that, for every
solution $x(t)$, $t\geq 0$, it holds that $x (t)\geq \varepsilon $ for
 all~$t\geq  \tau$.
\end{Lemma}
%%%%%%%%%%%%%%%%%%%%%%%%%%%

In other words, the conclusion is that after an arbitrarily short time
every~$x_i(t)$ is separated away from zero.

To prove Proposition~\ref{prop:bio}, note that
the Jacobian  of~\eqref{eq:bio} is
\be\label{eq:jac}
            J(x)= \begin{bmatrix}     -\alpha_1 &  0       & 0 &\dots&0&  g'(x_n)\\
                                              1 & -\alpha_2& 0 &\dots& 0 &0 \\
                                              0& 1 & -\alpha_3 &\dots& 0 &0 \\
                                              &&\vdots\\
                                              0& 0 & 0  &\dots& 1 &-\alpha_n                                     %%%%
                                                 \end{bmatrix}.
\ee
Thus,~$J(x)$
is a Metzler matrix,
so~\eqref{eq:bio} is a monotone dynamical system~\cite{hlsmith}.
 It is well-known~\cite[Ch.~3]{vid} that the induced
matrix measure corresponding to the~$L_1$ vector norm
  is
$
 \mu_1(A)=\max\{c_1(A),\ldots, c_n(A)\}
$,
 where \be\label{eq:cstac} c_j(A):=A_{jj}+\sum_{ \substack { 1\leq i \leq n\\  i \not = j} } |A_{ij}|  ,
 \ee
 %%%
  i.e.,   the sum of the entries in column~$j$ of~$A$,
 with non diagonal
elements   replaced by their absolute values. Of course, if~$A$ is Metzler then one  can take~$ A_{ij} $ instead  of~$|A_{ij}|$
  in~\eqref{eq:cstac}.
If~$P$ is an invertible matrix, and
$|\cdot|_{1,P} : \R^n \to \R_+$ is the vector norm defined by
$|z|_{1,P}:=|P  z|_1 $, then the
induced matrix measure is $
                    \mu_{1,P}(A):= \mu_1(PAP^{-1}).$
Let
\begin{align*}
D_\varepsilon :=  \diag \left(1,\alpha_1-\varepsilon , (\alpha_1-\varepsilon)(\alpha_2-\varepsilon ),\dots, \prod_{i=1}^{n-1}(\alpha_{i}-\varepsilon ) \right),
\end{align*}
 with~$\varepsilon>0$ sufficiently small.  Then
 %%%%%%%%%%%%%%%%%%%%
\begin{align*}
 D_\varepsilon & J(x) D_\varepsilon ^{-1}\\&= \begin{bmatrix}
   -\alpha_1 & 0 & 0  &\dots & 0 &  \frac{ g'(x_n) }{ \prod_{i=1}^{n-1}(\alpha_{i}-\varepsilon )  } \\
%%%%
  \alpha_1-\varepsilon &  -\alpha_2 &0& \dots &0 &0\\
  0 &   \alpha_2 -\varepsilon &0& \dots &0 &0\\
  &&\vdots\\
  0 &  0   &0& \dots &\alpha_{n-1}-\varepsilon  &-\alpha_n
%%%%%
    \end{bmatrix},
\end{align*}
so
\begin{align}\label{eq:2term}
\mu_{1,D_\varepsilon}(J(x))& =\max\{ -  \varepsilon  , \frac{ g'(x_n) }{\prod_{i=1}^{n-1}(\alpha_{i}-\varepsilon )  } -\alpha_n \}\nonumber \\
            &=\max\{   -  \varepsilon, \frac{ g'(x_n)-\alpha_n \prod_{i=1}^{n-1}(\alpha_{i}-\varepsilon )}
            {\prod_{i=1}^{n-1}(\alpha_{i}-\varepsilon )  }
             \}.
\end{align}
 Suppose that~$ {k-1}  < \alpha {k^2}$. Then for all~$x\in\R^n_+$,
\[
        g'(x_n)= \frac{k-1}{(k+x_n)^2} \leq \frac{k-1}{k^2}<\alpha .
\]
Combining this with~\eqref{eq:2term} implies that
there exists a sufficiently small~$\varepsilon>0$ such that
 $\mu_{1,D_\varepsilon}(J(x)) <-\varepsilon/2$
for all~$x\in \R^n_+$,  so the system is contractive on~$\R^n_+$.

Now assume that~${k-1} =\alpha {k^2}$.
By~\eqref{eq:jac},
\[
            \det(J(x))=(-1)^n (\alpha-g'(x_n)),
\]
so for every~$x \in \R^n_+$ with~$x_n=0$, we have~$\det(J(x))=(-1)^n (\alpha-g'(0))=0$.
This implies that the system is not contractive, with respect to any norm, on~$\R^n_+$.

We now use Theorem~\ref{thm:qcon} to
prove that~\eqref{eq:bio}  is {\sostshort}.
Note that since~$g'(u)=\frac{k-1}{(k+u)^2}$ and~$k>1$,
\[
            g(x_n)\geq g(0)=1/k,\quad \text{for all } x \in \R^n_+.
\]
For~$\zeta \in (0,1/2]$, let
\[
            \Omega_\zeta:=\{x \in \R^n_+: x\geq \zeta\}.
\]
It is straightforward to verify that~\eqref{eq:bio} satisfies   condition~(BR)
in  Lemma~\ref{lem:per}.
Hence, for every~$\tau>0$ there exists~$\varepsilon (\tau)>0$ such that
$x(t)\in \Omega_\varepsilon $ for all~$t\geq \tau$. Then
\[
        g'(x_n)= \frac{k-1}{(k+x_n)^2} \leq \frac{k-1}{(k+\varepsilon)^2}< \frac{k-1}{k^2} = \alpha .
\]
We already showed that this implies that there exists a~$\zeta>0$
and a norm~$|\cdot|_{1, D_\zeta }$
such that~\eqref{eq:bio} is contractive on~$\Omega_\varepsilon$ with respect to this norm.
When~$\zeta =0$,
\eqref{eq:2term} yields
$
\mu_{1,D_0}(J(x))
              =\max\{   0,  \frac{  g'(x_n)-\alpha   }{\alpha_{n-1} \dots \alpha_2\alpha_1 } \}
%%%%%
$,
  and since~$g'(x_n)  < \alpha$,
  $\mu_{1,D_0}(J(x)) \equiv 0$. Thus,~\eqref{eq:bio} is~NE with respect to~$|\cdot|_{1, D_0 }$.
  Summarizing, all the conditions in Theorem~\ref{thm:qcon} hold, and thus the system
  is~{\sostshort}. By Remark~\ref{rem:Global_Lip}, this implies~\soshort.~$\square$

%%%%%%%%%%%%%%%%%%%%%%%%%%%%%%%%%%%%%%%%
\noindent {\sl Proof of Lemma~\ref{lem:eqdef}.}\\
%%%%%%%%%%%%%%%%%%%%%%%%%%%%%%%%%%%%%%%%%%
If~\eqref{eq:fdyn}  is {\sostshort} then~\eqref{eq:qcontop} holds  for the particular case~$\varepsilon=\tau$ in Definition~\ref{def:qcont}. To prove the converse
implication, assume that~\eqref{eq:qcontop} holds. Pick~$\hat{\tau},\hat{\varepsilon}>0$. Let \be \label{eq:deftaumin} \tau:=\min\{\hat{\tau},\hat{\varepsilon}    \}, \ee and
let~$\ell=\ell(\tau)>0$. Pick~$t\geq t_1\geq 0$, and let~$t_2:=t+\hat{\tau}-\tau \geq t_1$. Then
%%%%%
 \begin{align*}
            |x(t_2 + {\tau},& t_1,a)-  x(t_2+\tau  ,t_1,b)|_v \nonumber \\
                    &\leq   (1+  {\tau}) \exp(-  ( t_2   -  t_1)  {\ell} ) |a-b|_v \\
                    &\leq (1+\hat{\varepsilon}) \exp(-  ( t   -  t_1)  {\ell} ) |a-b|_v,
%%%
 \end{align*}
where the last inequality follows from~\eqref{eq:deftaumin}. Thus,
%%%%%%%
 \begin{align*}
            |x(t +\hat{\tau},& t_1,a)-  x(t +\hat{\tau},t_1,b)|_v \nonumber \\
                    &\leq (1+\hat{\varepsilon}) \exp(-  ( t   -  t_1)  {\ell} ) |a-b|_v,
%%%
 \end{align*}
  and recalling that~$\hat{\tau},\hat{\varepsilon}>0$ were arbitrary,
 we conclude that Condition~2) in Lemma~\ref{lem:eqdef} implies
 \sostshort.

To prove that Condition~3) is equivalent to {\sostshort}, suppose that~\eqref{eq:qcontnew} holds. Then for any~$t_2\geq t_1$,
 \begin{align*}
                |x(t_2+\tau,t_1,a)-&x(t_2+\tau,t_1,b)|_v\\&\leq (1+\varepsilon) \exp(-  (t_2+\tau -t_1) \ell_1 ) |a-b|_v \\
                                                        &\leq (1+\varepsilon) \exp(-  (t_2  -t_1) \ell_1 ) |a-b| _v,
\end{align*}
so we have {\sostshort}.
%%%%%%%%%%%%
 Conversely, suppose that~\eqref{eq:fdyn}  is {\sostshort}.
 Pick any~$\tau,\varepsilon>0$. Then
 there exists~$\ell=\ell( {\tau}, {\varepsilon}/2  )>0$
 such that for any~$t \geq t_1+{\tau}$
 \begin{align*}
              |x(t ,t_1,a)-& x(t ,t_1,b)| _v\\
                                       & =  |x(t-{\tau}+{\tau},t_1,a)-  x(t-{\tau}+{\tau}  ,t_1,b)| _v  \\
              & \leq (1+ {\varepsilon}/2) \exp(-  (t-{\tau}  -t_1) \ell  ) |a-b|_v.
%%%             & =(1+2{\varepsilon})\exp(\tau \ell)  \exp(-  (t   -t_1) \ell  ) |a-b|.
%%%%%
 \end{align*}
 Thus, for any~$c \in(0,1)$
 \begin{align*}
              |x(t ,t_1,a)-& x(t ,t_1,b)| _v\\
             %% & \leq (1+ {\varepsilon}/2) \exp(-  (t-{\tau}  -t_1) c \ell  ) |a-b|_v\\
              & \leq (1+ {\varepsilon}/2) \exp( \tau c \ell ) \exp(-  (t   -t_1) c \ell  ) |a-b|_v.
%%%%%
 \end{align*}
Taking~$c>0$ sufficiently
 small such that~$(1+ {\varepsilon}/2) \exp( \tau c \ell ) \leq 1+\varepsilon$ implies that~\eqref{eq:qcontnew} holds for~$\ell_1:=c \ell$.
  This completes the proof
  that~\eqref{eq:qcontnew}
  is equivalent to {\sostshort}.~$\square$

%%%%%%%%%%%%%%%%%%%

%%%%%%%%%%%%%%%%%%%%%%%%%%%
\noindent {\sl Proof of Proposition~\ref{prop:sepimp}.}\\
%%%%%%%%%%%%%%%
Suppose that~\eqref{eq:fdyn} is  {\sostshort} with respect to some
norm~$|\cdot|_v$.
 Pick~$\varepsilon>0$.  Since the system is~{\sweshort}, there exists~$\tau_0=\tau_0(\varepsilon)>0$
such that
\[
                                        | x(t ,t_0,a)-x(t ,t_0,b) |_v
                                        \leq(1+ {\varepsilon}/{2})|a-b|_v, \
 \]
 for all $t \in [t_0,t_0+\tau_0]$.
Letting~$\ell_2:=\frac{1}{\tau_0}\log(\frac {1+\varepsilon}{1+(\varepsilon/2)} )$   yields
 \be\label{eq:nnsep}
                                        | x(t ,t_0,a)-x(t ,t_0,b) |_v
                                        \leq(1+\varepsilon) \exp(-  (t -t_0) \ell_2 ) |a-b|_v,
\ee
 for all $ t \in [t_0,t_0+\tau_0]$.
By item~\ref{item3} in Lemma~\ref{lem:eqdef} there exists~$\ell_1=\ell_1(\tau_0,\varepsilon)>0$
 such that
 \[
 %%%%%
            |x(t  ,t_0,a)-x(t  ,t_0,b)|_v  \leq   (1+\varepsilon) \exp(-  (t -t_0) \ell_1 ) |a-b|_v ,
 \]
 for all $t  \geq t_0+\tau_0$.
Combining this with~\eqref{eq:nnsep} yields
\[
 %%%%%
            |x(t  ,t_0,a)-x(t  ,t_0,b)|_v  \leq    (1+\varepsilon) \exp(-  (t -t_0) \ell  ) |a-b|_v ,
 \]
 for all $t  \geq t_0$,
            where~$\ell:=\min\{\ell_1,\ell_2\}>0$. This proves {\soshort}.~$\square$
%%%%%%%%%% end of proof

%%%%%%%%

%%%

 %%%%%%%%%%%%%%%%%%%%%%%%%%%%%%%%%%%%%%%%
\section{Conclusions}
%% %%%%%%%%%%%%%%%%%%%%%%%%%%%%%%%%%%%%%%%%%%%

Contraction theory has proved useful for studying numerous dynamical systems.
Contraction implies several desirable asymptotic  properties such as convergence to a unique attractor (if it exists)
 and entrainment to
periodic excitation. However, proving contraction is in many cases non-trivial.

In this paper,
we introduced  three   generalizations of contraction.
These are motivated by allowing contraction to take place
after an arbitrarily small transient
in  time and/or amplitude. We provided
conditions guaranteeing  that these forms of GC  hold, and
  demonstrated their usefulness
by  using them to analyze
 systems  that are
  not   contractive,  with respect to  any norm, yet are a~GC.

We note in passing that our original motivation for   generalizating
    contraction   was to
prove entrainment   in a model for translation-elongation
 called the \emph{ribosome flow model}~(RFM)~\cite{reuveni} (see also~\cite{RFM_stability,RFM_feedback,HRFM_steady_state,infi_HRFM}).
The state-variables~$x_i(t)$, $i=1,\dots,n$,
 in the RFM represent     occupancy levels on a coarse-grained model of the mRNA,
  normalized so that~$x_i(t)\in[0,1]$ for all~$t$.
  The
state-space of the  RFM is thus~$C^n:=[0,1]^n$.
It is straightforward to show, using
the  results presented here, that the~RFM is
 not contractive with respect to any norm on~$C^n$, yet
 is  {\stshort} on~$C^n$.

\subsubsection*{Acknowledgements}
%%%%%%%%%%%%%%%%%%%%%%%%%%%%%%%%%%%%%%%%%%
We thank Zvi Artstein for   helpful comments.

\bibliographystyle{IEEEtran}

\bibliography{cast_bib}

\begin{thebibliography}{10}
\providecommand{\url}[1]{#1}
\csname url@rmstyle\endcsname
\providecommand{\newblock}{\relax}
\providecommand{\bibinfo}[2]{#2}
\providecommand\BIBentrySTDinterwordspacing{\spaceskip=0pt\relax}
\providecommand\BIBentryALTinterwordstretchfactor{4}
\providecommand\BIBentryALTinterwordspacing{\spaceskip=\fontdimen2\font plus
\BIBentryALTinterwordstretchfactor\fontdimen3\font minus
  \fontdimen4\font\relax}
\providecommand\BIBforeignlanguage[2]{{%
\expandafter\ifx\csname l@#1\endcsname\relax
\typeout{** WARNING: IEEEtran.bst: No hyphenation pattern has been}%
\typeout{** loaded for the language `#1'. Using the pattern for}%
\typeout{** the default language instead.}%
\else
\language=\csname l@#1\endcsname
\fi
#2}}

\bibitem{LOHMILLER1998683}
W.~Lohmiller and J.-J.~E. Slotine, ``On contraction analysis for non-linear
  systems,'' \emph{Automatica}, vol.~34, pp. 683--696, 1998.

\bibitem{entrain2011}
G.~Russo, M.~di~Bernardo, and E.~D. Sontag, ``Global entrainment of
  transcriptional systems to periodic inputs,'' \emph{PLoS Computational
  Biology}, vol.~6, p. e1000739, 2010.

\bibitem{cont_mech}
W.~Lohmiller and J.-J.~E. Slotine, ``Control system design for mechanical
  systems using contraction theory,'' \emph{IEEE Trans.\ Automat.\ Control},
  vol.~45, pp. 984--989, 2000.

\bibitem{observer_posi_2011}
S.~Bonnabel, A.~Astolfi, and R.~Sepulchre, ``Contraction and observer design on
  cones,'' in \emph{Proc.\ 50th IEEE Conf. on Decision and Control and European
  Control Conference}, Orlando, Florida, 2011, pp. 7147--7151.

\bibitem{wang_slotine_2005}
W.~Wang and J.~J. Slotine, ``On partial contraction analysis for coupled
  nonlinear oscillators,'' \emph{Biol. Cybern.}, vol.~92, pp. 38--53, 2005.

\bibitem{partial_cont}
J.-J.~E. Slotine, ``Modular stability tools for distributed computation and
  control,'' \emph{Int. J. Adaptive Control and Signal Processing}, vol.~17,
  pp. 397--416, 2003.

\bibitem{russo_hier}
G.~Russo, M.~di~Bernardo, and E.~Sontag, ``A contraction approach to the
  hierarchical analysis and design of networked systems,'' \emph{IEEE Trans.\
  Automat.\ Control}, vol.~58, pp. 1328--1331, 2013.

\bibitem{angeli_inc}
D.~Angeli, ``A {L}yapunov approach to incremental stability properties,''
  \emph{IEEE Trans.\ Automat.\ Control}, vol.~47, pp. 410--421, 2002.

\bibitem{contra_sep}
\BIBentryALTinterwordspacing
F.~Forni and R.~Sepulchre, ``A differential {L}yapunov framework for
  contraction analysis,'' 2013. [Online]. Available:
  \url{http://arxiv.org/abs/1208.2943v2}
\BIBentrySTDinterwordspacing

\bibitem{Aminzare201331}
Z.~Aminzare and E.~D. Sontag, ``Logarithmic {L}ipschitz norms and
  diffusion-induced instability,'' \emph{Nonlinear Analysis: Theory, Methods \&
  Applications}, vol.~83, pp. 31--49, 2013.

\bibitem{slotine_cont_noise}
Q.-C. Pham, N.~Tabareau, and J.-J. Slotine, ``A contraction theory approach to
  stochastic incremental stability,'' \emph{IEEE Trans.\ Automat.\ Control},
  vol.~54, pp. 816--820, 2009.

\bibitem{soderling_survey}
G.~Soderlind, ``The logarithmic norm. {H}istory and modern theory,'' \emph{BIT
  Numerical Mathematics}, vol.~46, pp. 631--652, 2006.

\bibitem{cont_anc}
J.~Jouffroy, ``Some ancestors of contraction analysis,'' in \emph{Proc.\ 44th
  IEEE Conf. on Decision and Control}, Seville, Spain, 2005, pp. 5450--5455.

\bibitem{RufferWouwMueller:2013:Convergent-Systems-vs.-Incremental-Stabi:}
B.~S. R{\"u}ffer, N.~van~de Wouw, and M.~Mueller, ``Convergent systems vs.
  incremental stability,'' \emph{Systems Control Lett.}, vol.~62, pp. 277--285,
  2013.

\bibitem{RFM_entrain}
M.~Margaliot, E.~D. Sontag, and T.~Tuller, ``Entrainment to periodic initiation
  and transition rates in a computational model for gene translation,''
  \emph{PLoS ONE}, vol.~9, no.~5, p. e96039, 2014.

\bibitem{dorf-bullo}
F.~Dorfler and F.~Bullo, ``Synchronization and transient stability in power
  networks and nonuniform {K}uramoto oscillators,'' \emph{SIAM J.\ Control
  Optim.}, vol.~50, pp. 1616--1642, 2012.

\bibitem{hlsmith}
H.~L. Smith, \emph{Monotone Dynamical Systems: An Introduction to the Theory of
  Competitive and Cooperative Systems}, ser. Mathematical Surveys and
  Monographs.\hskip 1em plus 0.5em minus 0.4em\relax Providence, RI: Amer.
  Math. Soc., 1995, vol.~41.

\bibitem{sontag_textbook}
E.~D. Sontag, \emph{Mathematical Control Theory: Deterministic
  Finite-Dimensional Systems}, 2nd~ed., ser. Texts in Applied
  Mathematics.\hskip 1em plus 0.5em minus 0.4em\relax New York:
  Springer-Verlag, 1998, vol.~6.

\bibitem{mcs_angeli_2003}
D.~Angeli and E.~D. Sontag, ``Monotone control systems,'' \emph{IEEE Trans.\
  Automat.\ Control}, vol.~48, pp. 1684--1698, 2003.

\bibitem{vid}
M.~Vidyasagar, \emph{Nonlinear Systems Analysis}.\hskip 1em plus 0.5em minus
  0.4em\relax Englewood Cliffs, NJ: Prentice Hall, 1978.

\bibitem{reuveni}
S.~Reuveni, I.~Meilijson, M.~Kupiec, E.~Ruppin, and T.~Tuller, ``Genome-scale
  analysis of translation elongation with a ribosome flow model,'' \emph{PLoS
  Computational Biology}, vol.~7, p. e1002127, 2011.

\bibitem{RFM_stability}
M.~Margaliot and T.~Tuller, ``Stability analysis of the ribosome flow model,''
  \emph{IEEE/ACM Trans. Computational Biology and Bioinformatics}, vol.~9, pp.
  1545--1552, 2012.

\bibitem{RFM_feedback}
M.~Margaliot and T.~Tuller, ``Ribosome flow model with positive feedback,''
  \emph{J. Royal Society Interface}, vol.~10, p. 20130267, 2013.

\bibitem{HRFM_steady_state}
M.~Margaliot and T.~Tuller, ``On the steady-state distribution in the
  homogeneous ribosome flow model,'' \emph{IEEE/ACM Trans. Computational
  Biology and Bioinformatics}, vol.~9, pp. 1724--1736, 2012.

\bibitem{infi_HRFM}
Y.~Zarai, M.~Margaliot, and T.~Tuller, ``Explicit expression for the steady
  state translation rate in the infinite-dimensional homogeneous ribosome flow
  model,'' \emph{IEEE/ACM Trans. Computational Biology and Bioinformatics},
  vol.~10, pp. 1322--1328, 2013.

\end{thebibliography}

\end{document}